\documentclass{elsart}
\begin{document}

\begin{frontmatter}

\title{Remarks on Cantor's diagonalization proof of 1891}
\author[Slavica]{Slavica Vlahovic} 
and 
\author[Branislav]{Branislav Vlahovic\corauthref{cor}}
\ead{vlahovic@wpo.nccu.edu}
\corauth[cor]{Corresponding author.}

\address[Slavica]{Gunduliceva 2, Sisak, Croatia}
\address[Branislav]{North Carolina Central University, Durham, NC 27707, USA}





\begin{abstract}



 Remarks on the Cantor's nondenumerability proof of 1891 that the real 
 numbers are noncountable will be given. By the Cantor's diagonal procedure, 
 it is not possible to build numbers that are different from all numbers in a 
 general assumed denumerable sequence of all real numbers. The numbers created on 
 the diagonal of the assumed sequence are not different from the numbers in the 
 assumed denumerable sequence or they do not belong to that sequence.

\end{abstract}


\end{frontmatter}

\section{Introduction}

	The first proof that it is impossible to establish a one-to-one 
correspondence between the natural numbers $N$ and the real numbers $\Re$ is 
older than a century. In December of 1873 Cantor first proved non-denumerability of 
continuum and that first proof proceeded as follows\cite{1,2,3,4}: Find a closed 
interval $I_0$ that fails to contain $r_0$ then find a closed subinterval $I_1$ 
of $I_0$ such that $I_1$ misses $r_1$; continue in this manner, obtaining an 
infinite nested sequence of closed intervals, $I_0 \supseteq I_1 \supseteq I_2 
\supseteq ...$, that eventually excludes every one of the $r_n$; now let $d$ be 
a point lying in the intersection of all the Ia's; $d$ is a real number 
different from all of the $r_n$. 

This proof that no denumerable sequence of elements of an interval (a,b) can 
contain all elements of (a,b) often is overlooked in favor of the 1891 diagonal 
argument\cite{5}, when reference is made to Cantor's proving the 
nondenumerability of the continuum. Cantor himself repeated this proof with some 
modifications\cite{2,3,6,7,8,9,10,11,12,13,14} from 1874 to 1897, and today 
we have even more variations of this proof given by other authors. However, we 
have to note that they are in nuce similar; all of them include same modification 
of the Cantor's idea to derive a contradiction by defining in terms which cannot 
possibly be in the assumed denumerable sequence. So, in principle, all these proofs 
do not represent a significant change from Cantor's original idea and we can take them to be the 
same as the Cantor's proofs.

	For the reason of clarity, we will not discuss objections to these 
proofs that have been raised earlier\cite{15,16,17,18,19,20,21} or the 
legitimacy of these proofs from intuitionistic points of view \cite{22} 
and their nonconstructive parts, namely appeal to the Bolzano-Weierstrass  
theorem\cite{23} and inclusion of impredicative methods\cite{24}.  We will focus 
to show what is in principle wrong with the Cantor's  1891 proof and 
consequently all other similar proofs.

\section{Remarks on Cantor's 1891 diagonal proof of the nondenumerability of real numbers}

\vskip 0.3cm
{\bf Theorem 1}
\vskip 0.2cm

By the Cantor's diagonal procedure, it is not possible to build numbers that are 
different from all numbers in a general assumed denumerable sequence of all real 
numbers or created real numbers do not belong to the considered interval.  

\vskip 0.2cm
	{\bf Proof of the theorem 1}
	\vskip 0.2cm
	
	Cantor famous method of diagonalization is relaying upon only two elements, $m$ and $w$. With these he considered the collection A of elements $E = (x_1, x_2, ..., x_n, ...)$, where each $x_n$ was either $m$ or $w$. As example:
	
	$M = (m,m,m,m,...),$
	
	$W = (w,w,w,w,...),$
	
	$Emw = (m,w,m,w,...)$.
	
	Cantor then  asserted that the collection of all such elements $A$ was nondenumerable. 
	
	Let as repeat that proof by considering an open interval of numbers (M,W).  
	Cantor first produced a countable listing of elements $E_\nu$ in 
terms of the corresponding array (\ref{ar1}), where each $a_{\mu,\nu}$ was 
either m or w:

 \begin{equation}
 {\begin{array}{clcr}
     E_1=(a_{11},a_{12},...,a_{1\nu},...)\\
     E_2=(a_{21},a_{22},...,a_{2\nu},...)\\ 
     \vdots 
\;\;\;\;\;\;\;\;\;\;\;\;\;\;\;\;\;\;\;\;\;\;\;\;\;\;\;\;\;\;\;\;\;\;\; \\
     E_{\nu}=(a_{\mu1},a_{\mu2},...,a_{\mu\nu},...)\\
     \vdots 
\;\;\;\;\;\;\;\;\;\;\;\;\;\;\;\;\;\;\;\;\;\;\;\;\;\;\;\;\;\;\;\;\;\;\;
     \end{array}}
       \label{ar1}
 \end{equation}
       
       Then Cantor defined a new sequence $b_1, b_2,...,b_\nu,...$, where each $b_\nu$ 
was either m or w, determined so that $b_\nu\neq a_{\nu\nu}$. By formulating 
from this sequence of $b_\nu$ the element $E_0=(b_1,b_2,...,b_\nu,...)$, it 
followed that $E_0 \neq E_\nu$ for any value of the index $\nu$.  

However, this statements which appears so obvious, that whichever element $E_\nu$ one might choose to consider, there exists number $E_0$, which belongs to sequence (\ref{ar1}), and which is always different in $\nu^{th}$ coordinate, is not correct.  

By the Cantor, the number constructed on the diagonal must satisfy that $b_\nu\neq a_{\nu\nu}$. 
But the sequence (\ref{ar1}) might be arranged so that all $a_{\nu\nu} = m$. Therefore, in that case on the diagonal only one number might be created, which is $b_\nu = W$ However, number $W$ is not inside the interval $(M,W)$, so it is not required to be the part of the sequence (\ref{ar1}). It is obvious, that in this case Cantor can not establish contradiction, stating that there exists a number that should be part of the sequence (\ref{ar1}), but it is not listed in that sequence. This proves that Cantor's theorem is not correct.

\section{Conclusion}

It is impossible by the proposed diagonal procedure to build numbers that are not 
included in the assumed denumerable set and particularly it is not possible by this way 
to create an ascending hierarchy, in fact a limitless sequence of transfinite 
powers.

  	  	\end{document}